# Hopf (Bi-)Modules and Crossed Modules in Braided Monoidal Categories

Yuri Bespalov *    Bernhard Drabant †

September 1995

**Abstract**

Hopf (bi-)modules and crossed modules over a bialgebra $B$ in a braided monoidal category $\mathcal{C}$ are considered. The (braided) monoidal equivalence of both categories is proved provided $B$ is a Hopf algebra (with invertible antipode). Bialgebra projections and Hopf bimodule bialgebras over a Hopf algebra in $\mathcal{C}$ are found to be isomorphic categories. A generalization of the Majid-Radford criterion for a braided Hopf algebra to be a cross product is obtained as an application of these results.

## 1 Introduction

For bialgebras over a field $k$ the smash product and the smash coproduct are investigated extensively in the literature [Rad, Mol]. Let $H$ be a bialgebra, $B$ be an $H$-right module algebra and an $H$-right comodule coalgebra. If the smash product algebra structure and the smash coproduct coalgebra structure on $H \otimes B$ form a bialgebra then Radford [Rad] speaks of an admissible pair $(H, B)$. It is noted by Majid [Ma1] that in the work [Rad] the notion of a crossed module implicitly occurs as one of the conditions for a pair $(H, B)$ to be admissible. Crossed modules have been explicitely investigated in [Yet, RT] where the connection to knot theory and quantum groups was exhibited. Hopf bimodules are a special form of Hopf modules and appear as the basic notion in Woronowicz's approach to differential calculi on quantum groups [Wor] where they are called bicovariant bimodules. A connection of Hopf bimodules and crossed modules in terms of a choosen basis over $k$ has been found in [Wor]. For a symmetric monoidal category which admits (co-)equalizers, a coordinate free version of Woronowicz's result was found in [Sch]. The main theorem in [Sch] states the equivalence of the category of Hopf bimodules and the category of crossed modules over a particular Hopf algebra.

In the present paper we generalize this result to arbitrary braided monoidal (base) categories $\mathcal{C}$. The fundamental theorem of Hopf modules in braided categories [Lyu2] will be proven without the assumption of the existence of (co-)equalizers but with the help of the weaker condition that idempotents in $\mathcal{C}$ split. We recall the results of [B2] on crossed modules and define Hopf bimodules in the braided monoidal category $\mathcal{C}$. Then we prove in particular the (pre-)braided monoidal equivalence of the category of Hopf bimodules and the category of crossed modules over a Hopf algebra $H$ in $\mathcal{C}$. We find a connection of the bialgebra projections over a Hopf algebra $H$ (with isomorphic antipode) and the bialgebras in the category of $H$-Hopf bimodules. Both categories are isomorphic. One of the consequences of this theorem is the very natural description of (braided) admissible pairs $(H, B)$ in terms of the objects $B$ which are just the bialgebras in the category of $H$-crossed modules. Another application of the results of this paper is the construction of (bicovariant) differential calculi on braided Hopf algebras in an abelian, braided monoidal category. This will be published elsewhere.

*Supported by the International Science Foundation under grant U4J200
†Supported in part by a DFG research fellowship under grant Dr-288/1-1

The paper is organized as follows. In Section 2 we recall basic notations on braided monoidal categories. We give a short introduction to the graphical calculus [Ma4, Ma5, Yet, B2, Dra]. In Section 3 the definition of a Hopf module over a bialgebra $B$ in the braided monoidal category $\mathcal{C}$ is given and under the assumption that $\mathcal{C}$ admits split idempotents we introduce for any Hopf module over the Hopf algebra $H$ the idempotent $\Pi$ which is fundamental for what follows. The Structure Theorem of Hopf modules [Swe] will be proved, which states that the Hopf modules over $H$ are braided monoidal equivalent to the category $\mathcal{C}$ itself. As a sidestep we give similar results for two-fold Hopf modules over $H$ which have been investigated for the symmetric case in [RT, Sch]. In Section 4 we recall results of [B2] on crossed modules in braided monoidal categories. We give the notion of Hopf bimodules over a bialgebra $B$. In the main theorem of Section 4 we prove the (pre-)braided monoidal equivalence of Hopf bimodules and crossed modules over a Hopf algebra $H$ in $\mathcal{C}$. This is a braided equivalence if the antipode of $H$ is an isomorphism. Section 5 is devoted to the investigation of bialgebra projections over the Hopf algebra $H$. We prove that the category of bialgebra projections and the category of Hopf bimodule bialgebras over $H$ are isomorphic. As a consequence of this theorem the $H$-admissible objects, i.e. the objects $B$ in $\mathcal{C}$ which together with $H$ form a braided admissible pair $(H, B)$, are found to be the bialgebras in the category of crossed modules over $H$. In the appendix we explain what we mean by canonical splitting of idempotents. This is a technical assumption which provides functoriality of our constructions.

## 2 Preliminaries

### 2.1

Throughout the paper $\mathcal{C} := (\mathcal{C}, \otimes, \underline{1}, \alpha, \rho, \lambda, \Psi)$ is a *braided monoidal category* [FY, JS] where $\otimes$ is the tensor product (bifunctor), $\underline{1}$ is the unit object, $\alpha$ is the natural isomorphism which rules the associativity of the tensor product, $\rho$ and $\lambda$ are the natural isomorphisms for the right and left tensor multiplication with the unit object respectively, and $\Psi$ is the braiding. By Mac Lane's coherence theorem [Mac], $\mathcal{C}$ is equivalent to a strict monoidal category, i.e. a category where $\alpha$, $\rho$ and $\lambda$ are identity morphisms. This allows us to neglect the morphisms $\alpha$, $\rho$ and $\lambda$ in most of the calculations. We suppose that the reader is familiar with the notion of (co-)algebras, (co-)modules and bi-(co-)modules in monoidal categories [Swe, Ma4, Ma2, Ma3]. We assume (co-)associativity and the existence of a (co-)unit henceforth. In a braided monoidal category the tensor product of two (co-)algebras is again a (co-)algebra; the multiplication $m_{U \otimes V}$ and the unit $\eta_{U \otimes V}$ of two algebras $(U, m_U, \eta_U)$ and $(V, m_V, \eta_V)$ is given through

$$m_{U \otimes V} = (m_U \otimes m_V) \circ (\mathrm{id}_U \otimes \Psi_{U,V} \otimes \mathrm{id}_V), \quad \eta_{U \otimes V} = \eta_U \otimes \eta_V. \tag{1}$$

The coalgebra structure of the tensor product of two coalgebras is obtained in the *dual symmetric* manner, i.e. by reversing the order of the composition of morphisms, and by replacing the multiplication m by the comultiplication $\Delta$ and the unit $\eta$ by the counit $\varepsilon$. A bialgebra $(B, \mathrm{m}, \eta, \Delta, \varepsilon)$ in a braided monoidal category $\mathcal{C}$ is an algebra $(B, \mathrm{m}, \eta)$ and a coalgebra $(B, \Delta, \varepsilon)$ where $\Delta$ and $\varepsilon$ are algebra morphisms [Ma2, Ma3]. A Hopf algebra $(H, \mathrm{m}, \eta, \Delta, \varepsilon, S)$ in $\mathcal{C}$ is a bialgebra together with the antipode $S : H \to H$ such that

$$\mathrm{m} \circ (\mathrm{id}_H \otimes S) \circ \Delta = \mathrm{m} \circ (S \otimes \mathrm{id}_H) \circ \Delta = \eta \circ \varepsilon \,.$$

Every bialgebra $(B, \mathrm{m}, \eta, \Delta, \varepsilon)$ in $\mathcal{C}$ is a bi-(co-)module through the *regular action* $\mathrm{m}$ and the *regular coaction* $\Delta$. The tensor product of two (co-)modules over a bialgebra $B$ admits two kinds of (co-)actions. For instance the *diagonal action* of two right modules $(X, \mu_r^X)$ and $(Y, \mu_r^Y)$ is given by

$$\mu_{d,r}^{X \otimes Y} = (\mu_r^X \otimes \mu_r^Y) \circ (\mathrm{id}_X \otimes \Psi_{Y\,B} \otimes \mathrm{id}_B) \circ (\mathrm{id}_X \otimes \mathrm{id}_Y \otimes \Delta) \,, \tag{2}$$

and the *action induced by $Y$* is given through

$$\mu_{i,r}^{X \otimes Y} = \mathrm{id}_X \otimes \mu_r^Y \,. \tag{3}$$

(Dually) analogue all other types of diagonal (co-)actions and induced (co-)actions are defined. In the following also the (co-)adjoint (co-)action will be used; for completenes we recall the basic properties of it. Let $H$ be a Hopf algebra in $\mathcal{C}$ and let $(X, \mu_r, \mu_l)$ be an $H$-bimodule then $X$ becomes a right $H$-module through the right adjoint action

$$_{\mathrm{ad}}\triangleleft := \mu_l \circ (\mathrm{id}_H \otimes \mu_r) \circ (\Psi_{X\,H} \otimes \mathrm{id}_H) \circ (\mathrm{id}_X \otimes (S \otimes \mathrm{id}_H) \circ \Delta) \,. \tag{4}$$

Similarly the left adjoint action is defined. The coadjoint coactions are obtained in the dual symmetric manner. If $A$ is an algebra and $f : H \to A$ is an algebra morphism then the algebra $A$ becomes an $H$-bimodule $(A, \mu_r^f, \mu_l^f)$ via pullback along $f$, $\mu_l^f = \mathrm{m}_A \circ (f \otimes \mathrm{id}_A)$ and $\mu_r^f = \mathrm{m}_A \circ (\mathrm{id}_A \otimes f)$. The corresponding right adjoint action (induced by $f$) will be denoted by $_{\mathrm{ad}_f}\triangleleft$ and the resulting right $H$-module algebra by $A_f := (A, \mathrm{m}, \eta, {}_{\mathrm{ad}_f}\triangleleft)$.

## 2.2

In what follows we often use graphical notation for morphisms of monoidal categories [Ma4, Ma5, Yet, B2, Dra]. The graphics for (co-)multiplication, (co-)unit, antipode, left and right (co-)action, and braiding is given in Figure 1.

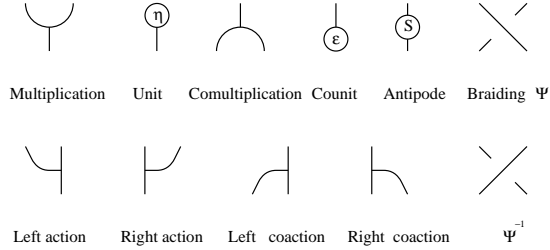

Multiplication   Unit   Comultiplication  Counit   Antipode   Braiding $\Psi$     (Fig. 1)

Left action   Right action   Left coaction   Right coaction   $\Psi^{-1}$

As an example the graphical representation of the associativity, the counit property, the bialgebra axiom for the multiplicativity of the comultiplication and the antipode property of an algebra, a coalgebra, a bialgebra, and a Hopf algebra respectively, are listed in Figure 2.

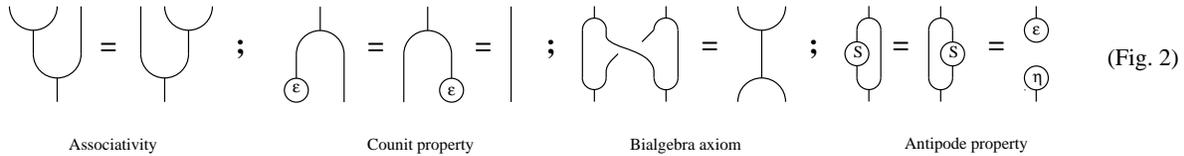

Associativity          Counit property          Bialgebra axiom          Antipode property      (Fig. 2)

**Remark 2.2.1** *If $(\mathcal{C}, \otimes, \underline{1}, \Psi)$ is a braided monoidal category then $\overline{\mathcal{C}} = (\mathcal{C}, \otimes, \underline{1}, \overline{\Psi})$ is braided monoidal where the mirror-reversed braiding $\overline{\Psi}_{X,Y} := \Psi_{Y,X}{}^{-1}$ is used. For a Hopf algebra $H$ in $\mathcal{C}$ we denote by $H^{\mathrm{op}}$ (resp. $H_{\mathrm{op}}$) the same coalgebra (resp. algebra) with the opposite multiplication $\overline{m}$ (resp. opposite comultiplication $\overline{\Delta}$) defined through*

$$\overline{m} := m \circ \Psi_{H\,H}{}^{-1} \qquad \left(\text{resp.} \qquad \overline{\Delta} := \Psi_{H\,H}{}^{-1} \circ \Delta\right). \tag{5}$$

*The opposite (co-)multiplication in (5) may not be confused with the opposite (co-)multiplication defined in [Ma4].*

**Lemma 2.2.1** *[Ma2] $H^{\mathrm{op}}$ and $H_{\mathrm{op}}$ are Hopf algebras in $\overline{\mathcal{C}}$ with antipode $S^{-1}$.*

# 3 Hopf modules

From [Lyu2] it is known that the Structure Theorem of Hopf modules [Swe] also holds for Hopf algebras in braided monoidal categories if there exist (co-)equalizers. We show in this section that to prove the Structure Theorem in a braided monoidal category $\mathcal{C}$ it suffices to assume that $\mathcal{C}$ admits split idempotents.

## 3.1

In a category $\mathcal{D}$ the idempotent $\mathrm{e} = \mathrm{e}^2 : X \to X$ is called split in $\mathcal{D}$ if there exists an object $X_{\mathrm{e}}$ and morphisms $\mathrm{i}_{\mathrm{e}} : X_{\mathrm{e}} \to X$ and $\mathrm{p}_{\mathrm{e}} : X \to X_{\mathrm{e}}$ such that $\mathrm{e} = \mathrm{i}_{\mathrm{e}} \circ \mathrm{p}_{\mathrm{e}}$ and $\mathrm{id}_{\mathrm{e}} = \mathrm{p}_{\mathrm{e}} \circ \mathrm{i}_{\mathrm{e}}$. If any idempotent e is split then $\mathcal{D}$ is said to admit split idempotents.

For a given category $\mathcal{C}$ there exists a universal category $\widehat{\mathcal{C}}$ which admits split idempotents (and is called *Karoubi enveloping category of* $\mathcal{C}$) and a full embedding $\mathcal{C} \xrightarrow{i} \widehat{\mathcal{C}}$ such that for any category $\mathcal{D}$ which split idempotents every functor $F : \mathcal{C} \to \mathcal{D}$ factorizes uniquely over $i$. $\widehat{\mathcal{C}}$ can be realized as the category with objects $X_e = (X, e)$, where $X$ is an object in $\mathcal{C}$ and $e : X \to X$ is idempotent, $e^2 = e$. The morphisms in $\widehat{\mathcal{C}}$ are defined by $\widehat{\mathcal{C}}(X_e, Y_f) := \{t \in \mathcal{C}(X, Y) \,|\, fte = t\}$. The full embedding $i$ is given through $i(X) = X_{\mathrm{id}_X}$, $i(f) = f$. Further information on the Karoubi enveloping category can be found in the appendix. It is noted in [Lyu1] that for $\mathcal{C}$ being a (braided) monoidal category the category $\widehat{\mathcal{C}}$ can be equipped with a (braided) monoidal structure.

$$\begin{aligned}
\underline{1} &:= (\underline{1}, \mathrm{id}_{\underline{1}}), \\
X_e \otimes Y_f &:= (X \otimes Y)_{e \otimes f}, \\
\alpha_{X_e, Y_f, Z_g} &:= ((e \otimes f) \otimes g) \circ \alpha_{X,Y,Z}, \\
\lambda_{X_e} &:= e \circ \lambda_X, \quad \rho_{X_e} := e \circ \rho_X, \\
\Psi_{X_e, Y_f} &:= (f \otimes e) \circ \Psi_{X,Y}.
\end{aligned}$$

In this case $i$ is a (braided) monoidal functor.

In a monoidal category $\mathcal{D}$ which admits split idempotents the category of (co-)modules over a (co-)algebra admits split idempotents. In anticipation of Section 4 similar facts hold for the categories of crossed modules and of Hopf (bi-)modules over a bialgebra in a braided monoidal category which admits split idempotents.

From now on we assume that the category $\mathcal{C}$ is braided monoidal and admits split idempotents. Moreover, we assume that for any idempotent $\mathrm{e} = \mathrm{e}^2 : X \to X$ an object $X_{\mathrm{e}}$ and morphisms $\mathrm{i}_{\mathrm{e}} : X_{\mathrm{e}} \to X$, $\mathrm{p}_{\mathrm{e}} : X \to X_{\mathrm{e}}$ are chosen ia a certain canonical way as explained in the appendix. This provides functoriality of our constructions.

## 3.2

**Definition 3.2.1** *A left Hopf module (resp. right-left Hopf module) $X$ over a bialgebra $B$ in $\mathcal{C}$ is a left (resp. right) $B$-module and a left $B$-comodule such that the action is a comodule morphism, where the regular (co-)action on $B$ and the diagonal tensor (co-)action on the tensor product of $X$ and $B$ is used. I.e. besides the (co-)module relations for $X$ a "polarized" version of the bialgebra axiom holds:*

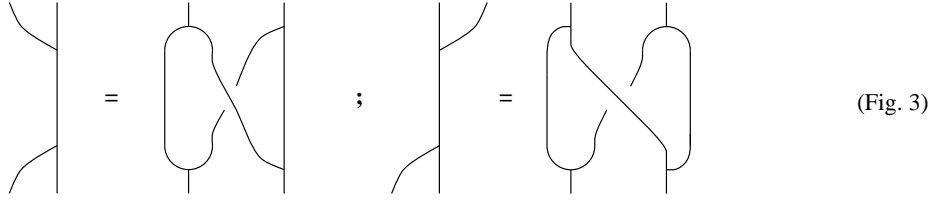

(Fig. 3)

*Left (resp. right-left) Hopf modules together with the left-(resp. right-)$B$-module-left-$B$-comodule morphisms form the category of Hopf modules which will be denoted by ${}_B^B\mathcal{C}$ (resp. ${}^B\mathcal{C}_B$). Similar definitions yield all other combinations of Hopf modules.*

It follows directly from the definition that any bialgebra $B$ is a Hopf module over itself with the regular action and the regular coaction.

**Remark 3.2.1** *A right-left $B$-Hopf module over a bialgebra $B$ in $\mathcal{C}$ is a left $B^{\mathrm{op}}$-Hopf module in $\overline{\mathcal{C}}$ where the opposite action is defined by the "polarized" form of the first equation in (5). Similarly one obtains from a left Hopf module all other versions of Hopf modules and therefore the results about left Hopf modules can be easily reformulated for any other type of Hopf modules. Hence we will restrict to left Hopf modules in the following.*

In the next proposition the idempotent $\Pi$ will be introduced which will be used later to construct a tensor product for Hopf (bi-)modules.

**Proposition 3.2.1** *Let $H$ be a Hopf algebra in $\mathcal{C}$ and $(X, \mu_l, \nu_l)$ be a left $H$-Hopf module. Then it holds:*

1. *The morphism ${}_X\Pi : X \to X$ defined through*

$$ {}_X\Pi := \mu_l \circ (S \otimes \mathrm{id}_X) \circ \nu_l \tag{6} $$

   *is an idempotent in $\mathrm{End}_{\mathcal{C}}(X)$.*

2. *Let ${}_H X \underset{{}_X\mathrm{p}}{\overset{{}_X\mathrm{i}}{\rightleftarrows}} X$ be the morphisms which split the idempotent ${}_X\Pi$, i.e. ${}_X\mathrm{i} \circ {}_X\mathrm{p} = {}_X\Pi$ and ${}_X\mathrm{p} \circ {}_X\mathrm{i} = \mathrm{id}_{{}_H X}$. Then*

$$ {}_H X \xrightarrow{{}_X\mathrm{i}} X \underset{\eta \otimes \mathrm{id}_X}{\overset{\nu_l}{\rightrightarrows}} H \otimes X \qquad \text{and} \qquad H \otimes X \underset{\varepsilon \otimes \mathrm{id}_X}{\overset{\mu_l}{\rightrightarrows}} X \xrightarrow{{}_X\mathrm{p}} {}_H X \tag{7} $$

   *are equalizer and coequalizer respectively. Hence ${}_H X$ is at the same time an object of invariants and coinvariants of $X$.*

**Proof.** Figure 4 shows that $\nu_l \circ {}_X\Pi = \eta \otimes {}_X\Pi$. In a dual way one obtains ${}_X\Pi \circ \mu_l = \varepsilon \otimes {}_X\Pi$. This immediately implies that ${}_X\Pi$ is an idempotent.

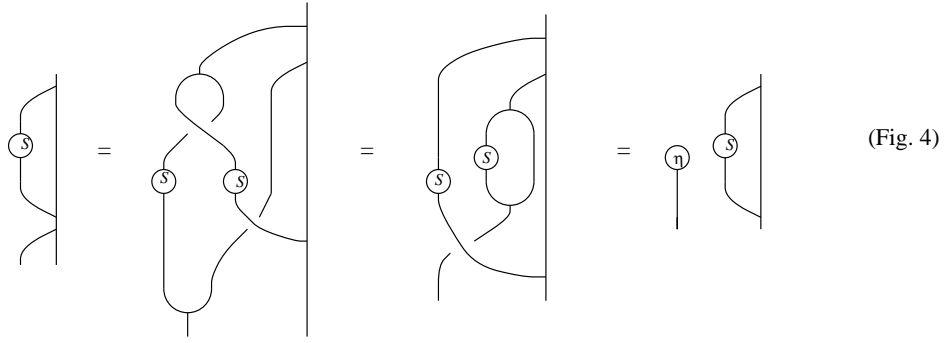

(Fig. 4)

Let $f: X' \to X$ be any morphism such that $\nu_l \circ f = \eta \otimes f: X' \to H \otimes X$. Then $f = {}_X\Pi \circ f = {}_Xi \circ ({}_Xp \circ f)$. Since ${}_Xi$ is monomorphic this proves the universality of ${}_HX \xrightarrow{{}_Xi} X$ as equalizer. In an analogous manner the universality of ${}_Xp$ as coequalizer will be proved. □

## 3.3

Before we investigate the relation between ${}_H^H\mathcal{C}$ and $\mathcal{C}$ we provide some useful results.

**Lemma 3.3.1** *Let $B$ be a bialgebra in $\mathcal{C}$, $X$ be a left $B$-Hopf module and $Y$ be an object in $\mathcal{C}$. Then $X \otimes Y$ is a $B$-Hopf module with action $\mu_l^{X \otimes Y} = \mu_l^X \otimes \mathrm{id}_Y$ and coaction $\nu_l^{X \otimes Y} = \nu_l^X \otimes \mathrm{id}_Y$ induced by $X$.*

**Proposition 3.3.2** *Let $B$ be a bialgebra in $\mathcal{C}$. Then the assignment $B \ltimes (-): \mathcal{C} \to {}_B^B\mathcal{C}$ with $B \ltimes (X) = (B \otimes X, \mu_l^{B \otimes X}, \nu_l^{B \otimes X})$, $\forall X \in \mathrm{Ob}(\mathcal{C})$ according to Lemma 3.3.1 and $B \ltimes (g) = \mathrm{id}_B \otimes g$ for morphisms in $\mathcal{C}$, defines a full inclusion functor.*

**Proof.** From Lemma 3.3.1 it is clear that $B \ltimes (X)$ is $B$-Hopf module. Obviously $B \ltimes (g)$ is Hopf module morphism, and the composition and the identity are preserved. Any Hopf module morphism $f: B \ltimes X \to B \ltimes Y$ has the form $f = \mathrm{id}_B \otimes f_B$ where $f_B := (\epsilon \otimes \mathrm{id}_Y) \circ f \circ (\eta \otimes \mathrm{id}_X)$:

$$\begin{aligned}
f &= f \circ (\mu \otimes \mathrm{id}_X) \circ (\mathrm{id}_B \otimes \eta \otimes \mathrm{id}_X) \\
&= (\mu \otimes \mathrm{id}_Y) \circ (\mathrm{id}_B \otimes f) \circ (\mathrm{id}_B \otimes \eta \otimes \mathrm{id}_X) \\
&= (\mu \otimes \epsilon \otimes \mathrm{id}_Y) \circ (\mathrm{id}_B \otimes \Delta \otimes \mathrm{id}_Y) \circ (\mathrm{id}_B \otimes f) \circ (\mathrm{id}_B \otimes \eta \otimes \mathrm{id}_X) \\
&= (\mu \otimes \epsilon \otimes \mathrm{id}_Y) \circ (\mathrm{id}_{B \otimes B} \otimes f) \circ (\mathrm{id}_B \otimes \Delta \otimes \mathrm{id}_X) \circ (\mathrm{id}_B \otimes \eta \otimes \mathrm{id}_X) \\
&= \mathrm{id}_B \otimes [(\epsilon \otimes \mathrm{id}_Y) \circ f \circ (\eta \otimes \mathrm{id}_X)]
\end{aligned}$$

The second (resp. the fourth) identity uses the fact that $f$ is a module (resp. comodule) morphism. Since $B$ is bialgebra it holds $B \ltimes f = B \ltimes g \iff f = g$. □

**Lemma 3.3.3** *Let $H$ be a Hopf algebra in $\mathcal{C}$ and $X$ be an $H$-Hopf module. According to Proposition 3.2.1 and Proposition 3.3.2 the tensor product $H \otimes {}_HX$ is an $H$-Hopf module. Then the morphisms*

$$\begin{aligned}
{}_X\mu &= \mu_l \circ (\mathrm{id}_H \otimes {}_Xi): H \otimes {}_HX \to X \\
{}_X\nu &= (\mathrm{id}_H \otimes {}_Xp) \circ \nu_l: X \to H \otimes {}_HX
\end{aligned} \qquad (8)$$

*are mutually inverse Hopf module morphisms.*

**Proof.** The proof of ${}_X\mu \circ {}_X\nu = \mathrm{id}_X$ is rather simple. Using the definition of a Hopf module and the properties of the morphisms ${}_Xi$ and ${}_Xp$ according to (7) one obtains

$$\begin{aligned}
{}_X\nu \circ {}_X\mu &= (\mathrm{id}_H \otimes {}_Xp) \circ (m \otimes \mu_\ell) \circ (\mathrm{id}_H \otimes \Psi_{H,H} \otimes \mathrm{id}_{{}_HX}) \circ (\Delta \otimes \nu_\ell) \circ (\mathrm{id}_H \otimes {}_Xi) \\
&= (m \circ (\mathrm{id}_H \otimes (\eta \circ \varepsilon))) \circ \Delta \otimes {}_Xp \circ {}_Xi = \mathrm{id}_A \otimes \mathrm{id}_{{}_HX}
\end{aligned}$$

The Hopf module properties of ${}_X\mu$ and ${}_X\nu$ are obvious. □

**Lemma 3.3.4** *Let $f : X \to Y$ be a Hopf module morphism over the Hopf algebra $H$. Then $f \circ {}_X\Pi = {}_Y\Pi \circ f$.*

Using these results yields

**Proposition 3.3.5** *Let $H$ be a Hopf algebra in $\mathcal{C}$, then the assignment ${}_H(-) : {}_H^H\mathcal{C} \to \mathcal{C}$ which is given through ${}_H(X) := {}_HX$ for an object $X$, and through ${}_H(f) = {}_Y\mathrm{p} \circ f \circ {}_X\mathrm{i}$ for a Hopf module morphism $f : X \to Y$ defines a functor.*

**Remark 3.3.1** 1. *Completely analogous results are obtained for the category $\mathcal{C}_H^H$, e.g. the idempotent $\Pi_X$, the functors $(-) \rtimes H$ and $(-)_H$, and the corresponding statements are obtained in a mirror reversed manner.*
2. *For a Hopf algebra $H$ it holds $\eta \circ \varepsilon = \mathrm{m} \circ (S \otimes \mathrm{id}_H) \circ \Delta = \mathrm{m} \circ (\mathrm{id}_H \otimes S) \circ \Delta = \Pi_H = {}_H\Pi$ and $\varepsilon \circ \eta = \mathrm{id}_{\underline{1}}$, hence ${}_HH = H_H = \underline{1}$.*

## 3.4

Let $A$ be an algebra in the category $\mathcal{C}$ then *the tensor product over $A$* of a right $A$-module $U$ and a left $A$-module $V$ is defined (up to isomorphism) as the coequalizer $\lambda_{U,V}^H : U \otimes V \to U \underset{A}{\otimes} V$ of the morphisms $\mu_r^U \otimes \mathrm{id}_V$ and $\mathrm{id}_U \otimes \mu_l^V$ if it exists. Dually to this construction is the notion of *the cotensor product over $C$* of a right $C$-comodule $P$ and a left $C$-comodule $Q$, denoted by $\rho_{P,Q}^C : P \underset{C}{\square} Q \to P \square Q$ if it exists, where $C$ is a coalgebra in $\mathcal{C}$.

The next proposition provides existence criteria for a (co-)tensor product over a Hopf algebra in $\mathcal{C}$.

**Proposition 3.4.1** *Let $H$ be a Hopf algebra, $M$ be a left $H$-Hopf module, $N$ be a right $H$-module and $P$ be a right $H$-comodule in $\mathcal{C}$. Then the tensor product of $N$ and $M$ over $H$ exists and is given through*

$$N \otimes M \xrightarrow{\lambda_{N,M}^H} N \underset{H}{\otimes} M \cong N \otimes {}_HM$$
$$\lambda_{N,M}^H = (\mu_r^N \otimes \mathrm{id}_{{}_HM}) \circ (\mathrm{id}_N \otimes {}_M\nu). \tag{9}$$

*The cotensor product of $P$ and $M$ over $H$ exists and is given by*

$$P \otimes {}_HM \cong P \underset{H}{\square} M \xrightarrow{\rho_{P,M}^H} P \otimes M$$
$$\rho_{P,M}^H = (\mathrm{id}_P \otimes {}_M\mu) \circ (\nu_r^P \otimes \mathrm{id}_{{}_HM}) \tag{10}$$

*where ${}_M\nu$ and ${}_M\mu$ according to (8) are used. If in addition $N$ is right $H$-Hopf module then $N \underset{H}{\otimes} M$ and $N \underset{H}{\square} M$ coincide.*

**Proof.** Using the results of Lemma 3.3.3 and the unit property of the action $\mu_r^N$, it is easily verified that $\lambda_{N,M}^H$ is epimorphic. To prove that $\lambda_{N,M}^H$ equalizes $\mu_r^N \otimes \mathrm{id}_M$ and $\mathrm{id}_N \otimes \mu_l^M$ one proceeds as follows.

$$\begin{aligned}
\lambda_{N,M}^H \circ (\mathrm{id}_N \otimes \mu_l^M) &= (\mu_r^N \otimes \mathrm{id}_{{}_HM}) \circ (\mathrm{id}_N \otimes {}_M\nu) \circ (\mathrm{id}_N \otimes \mu_l^M) \\
&= (\mu_r^N \circ (\mathrm{id}_N \otimes \mathrm{m}) \otimes \mathrm{id}_{{}_HM}) \circ (\mathrm{id}_{N \otimes H} \otimes {}_M\nu) \\
&= \lambda_{N,M}^H \circ (\mu_r^N \otimes \mathrm{id}_M)
\end{aligned} \tag{11}$$

The second equality in (11) uses the fact that ${}_M\nu$ is module morphism and the third equation can be obtained with the help of the associativity of the action $\mu_r^N$

Now suppose that $f$ is a morphism which also coequalizes $\mu_r^N \otimes \mathrm{id}_M$ and $\mathrm{id}_N \otimes \mu_l^M$. Then $f$ factorizes over $\lambda_{N,M}^H$.

$$\begin{aligned} f &= f \circ (\mathrm{id}_N \otimes {}_M\mu \circ {}_N\nu) \\ &= f \circ (\mu_r^N \otimes {}_M\mathrm{i}) \circ (\mathrm{id}_N \otimes {}_M\nu) \\ &= f \circ (\mathrm{id}_N \otimes {}_M\mathrm{i}) \circ \lambda_{N,M}^H \end{aligned}$$

where in the first equation Lemma 3.3.3 is used. The second equation holds by assumption for $f$. In the dual symmetric way is is proved that $\rho_{N,M}^H$ is the corresponding equalizer. $\square$

**Remark 3.4.1** *For a left $H$-Hopf module $M$ and a right $H$-Hopf module $N$ the composition morphism $\phi_{N,M}^H := \rho_{N,M}^H \circ \lambda_{N,M}^H$ equals to*

$$(\mu_r^N \otimes \mu_\ell^M) \circ (\mathrm{id}_N \otimes \Psi_{H\,H} \otimes \mathrm{id}_M) \circ (\nu_r^N \otimes \nu_\ell^M) \tag{12}$$

## 3.5

In case $H$ is a Hopf algebra in $\mathcal{C}$ it is possible to define a braided monoidal structure on the category ${}_H^H\mathcal{C}$.

**Theorem 3.5.1** *Let $H$ be a Hopf algebra in $\mathcal{C}$. Then $\otimes_H := \otimes \circ (\mathrm{id} \times {}_H(-)) : {}_H^H\mathcal{C} \times {}_H^H\mathcal{C} \to {}_H^H\mathcal{C}$ is a bifunctor. Through*

$$\,^{^H_H\mathcal{C}}\!\alpha_{X\,Y\,Z} := \alpha_{X\,{}_HY\,{}_HZ} \qquad \,^{^H_H\mathcal{C}}\!\rho_X := \rho_X \qquad \,^{^H_H\mathcal{C}}\!\lambda_X := {}_X\mu \quad \forall\, X, Y, Z \in {}_H^H\mathcal{C} \tag{13}$$

*the category $({}_H^H\mathcal{C}, H, \otimes_H, \,^{^H_H\mathcal{C}}\!\alpha, \,^{^H_H\mathcal{C}}\!\rho, \,^{^H_H\mathcal{C}}\!\lambda)$ is monoidal. ${}_H^H\mathcal{C}$ is a braided category by the braiding given through*

$$\,^{^H_H\mathcal{C}}\!\Psi_{X\,Y} = (\mu_l^Y \otimes {}_X\mathrm{p}) \circ (\mathrm{id}_H \otimes \Psi_{X\,Y}) \circ (\nu_l^X \otimes {}_Y\mathrm{i}) \tag{14}$$

*on objects $X, Y \in \mathrm{Ob}({}_H^H\mathcal{C})$.*

**Proof.** Lemma 3.3.1 and Proposition 3.3.5 immediately imply that $\otimes_H$ is well-defined. By Lemma 3.3.1 it is clear that

$$_{(X\otimes_H Y)}\mathrm{p} = {}_X\mathrm{p} \otimes \mathrm{id}_{{}_HY} \quad \text{and} \quad _{(X\otimes_H Y)}\mathrm{i} = {}_X\mathrm{i} \otimes \mathrm{id}_{{}_HY} \quad \forall\, X, Y \in \mathrm{Ob}({}_H^H\mathcal{C}) . \tag{15}$$

Using in particular (15) and Remark 3.3.1 a straightforward calculation shows that $\,^{^H_H\mathcal{C}}\!\alpha, \,^{^H_H\mathcal{C}}\!\rho$ and $\,^{^H_H\mathcal{C}}\!\lambda$ are the appropriate natural isomorphisms such that $({}_H^H\mathcal{C}, H, \otimes_H, \,^{^H_H\mathcal{C}}\!\alpha, \,^{^H_H\mathcal{C}}\!\rho, \,^{^H_H\mathcal{C}}\!\lambda)$ is a monoidal category. The morphism $\,^{^H_H\mathcal{C}}\!\Psi_{X\,Y}$ is isomorphic in $\mathcal{C}$. Its inverse is given by

$$(\,^{^H_H\mathcal{C}}\!\Psi_{X\,Y})^{-1} = (\mu_l^X \otimes {}_Y\mathrm{p}) \circ (\mathrm{id}_H \otimes (\Psi_{X\,Y})^{-1}) \circ (\nu_l^Y \otimes {}_X\mathrm{i}) . \tag{16}$$

This is easily checked with the help of the Hopf module property and Proposition 3.2.1.2. Exploiting Lemma 3.3.4 it follows that

$$\,^{^H_H\mathcal{C}}\!\Psi_{X'\,Y'} \circ (f \otimes_H g) = (g \otimes_H f) \circ \,^{^H_H\mathcal{C}}\!\Psi_{X\,Y} \tag{17}$$

where $f : X \to X'$ and $g : Y \to Y'$ are Hopf module morphisms. Lemma 3.3.1, the Hopf module axiom and Proposition 3.2.1 imply that $\,^{^H_H\mathcal{C}}\!\Psi_{X\,Y}$ is a left module morphism and dually analogue one shows that $\,^{^H_H\mathcal{C}}\!\Psi_{X\,Y}$ is a left comodule morphism. The defining relations for a braiding hold.

$$\begin{aligned} \,^{^H_H\mathcal{C}}\!\Psi_{X\,(Y\otimes_H Z)} &= (\mathrm{id}_Y \otimes_H \,^{^H_H\mathcal{C}}\!\Psi_{X\,Z}) \circ (\,^{^H_H\mathcal{C}}\!\Psi_{X\,Y} \otimes_H \mathrm{id}_Z) \\ \,^{^H_H\mathcal{C}}\!\Psi_{(X\otimes_H Y)\,Z} &= (\,^{^H_H\mathcal{C}}\!\Psi_{X\,Z} \otimes_H \mathrm{id}_Y) \circ (\mathrm{id}_X \otimes_H \,^{^H_H\mathcal{C}}\!\Psi_{Y\,Z}) . \end{aligned}$$

In fact using Lemma 3.3.1 and $_H\bigl(^{^H_H C}\Psi_{X\,Y}\bigr) = \Psi_{_HX\,_HY}$ yields the desired identity. □

These results provide us with the necessary techniques to proof the following theorem.

**Theorem 3.5.2** (Structure Theorem) *Let $H$ be a Hopf algebra in $\mathcal{C}$. Then the categories $\mathcal{C}$ and $^H_H\mathcal{C}$ are braided monoidal equivalent. The equivalence is given by $\mathcal{C} \underset{_H(-)}{\overset{H \ltimes (-)}{\rightleftarrows}} {}^H_H\mathcal{C}$.*

**Proof.** It is not difficult to prove that the natural isomorphism $\lambda_X : \underline{1} \otimes X \to X$ can be written as $\lambda : {}_H(-) \circ (H \ltimes (-)) \overset{\bullet}{=} \mathrm{id}_{\mathcal{C}}$. Similarly it can be checked that $_X\mu : H \otimes_H X \to X$ according to Lemma 3.3.3 is a natural isomorphism. Hence $\mathcal{C}$ and $^H_H\mathcal{C}$ are equivalent categories. It is a straightforward calculation to check that the natural isomorphims

$$\xi : \otimes_H \circ (H \ltimes (-), H \ltimes (-)) \overset{\bullet}{=} (H \ltimes (-)) \circ \otimes \tag{18}$$
$$\xi_{X\,Y} = \alpha^{-1}_{H\,X\,Y} \circ (\mathrm{id}_{H \ltimes X} \otimes \lambda_Y)$$

and

$$\mathrm{id} : {}_H(-) \circ \otimes_H \overset{\bullet}{=} \otimes \circ ({}_H(-), {}_H(-)) \tag{19}$$

are compatible (in the usual sense) with the unit objects, the associativity of the tensor products, the isomorphisms which rule the left and right unit tensor multiplication, and with the braidings. Thus $\mathcal{C} \cong {}^H_H\mathcal{C}$ are braided monoidal equivalent categories. □

## 3.6 Two-fold Hopf modules

Very similar to the previous considerations on Hopf modules are the results on two-fold Hopf modules. A generalization to Hopf bimodules, however, needs a little more care since one has to take into account the compatibility between right module and right comodule structure, but the results derived up to now induce the monoidal properties of Hopf bimodules for instance. This will be performed in detail in Section 4.

**Definition 3.6.1** *Let $B$ be a bialgebra in $\mathcal{C}$. A two-fold Hopf module $X = (X, \mu_l, \mu_r, \nu_l)$ is an object which is a $B$-bimodule in the category of left $B$-comodules, or in the language of Hopf modules: $X \in \mathrm{Ob}(^B_B\mathcal{C})$ and $X \in \mathrm{Ob}(^B\mathcal{C}_B)$.*
*Two-fold Hopf modules together with the $B$-bimodule-left $B$-comodule morphisms form the category of two-fold Hopf modules $^B_B\mathcal{C}_B$. Similarly all other types of two-fold Hopf modules are defined.*

**Proposition 3.6.1** *Let $B$ be a bialgebra in $\mathcal{C}$. Then for any $X \in {}^B_B\mathcal{C}_B$ and any $Y \in \mathcal{C}_B$ the object $X \otimes Y$ equipped with the left action and the left coaction induced by $X$ according to (3) and with the right diagonal action according to (2) is a two-fold Hopf module in $^B_B\mathcal{C}_B$. Like in Theorem 3.5.1 it is verified that $^B_B\mathcal{C}_B$ is monoidal with tensor product $\otimes_B$.*

Considering the regular object $B \in {}^B_B\mathcal{C}_B$ the previous construction extends to the functor $B \ltimes (-) : \mathcal{C}_B \to {}^B_B\mathcal{C}_B$ between monoidal categories. Using the adjoint action according to (4) yields

**Proposition 3.6.2** *Let $H$ be a Hopf algebra in $\mathcal{C}$. The following identities hold for a two-fold Hopf module $X \in \mathrm{Ob}(^H_H\mathcal{C}_H)$:*

$$_X\Pi \circ (_{\mathrm{ad}}\triangleleft) = {}_X\Pi \circ \mu_r = (_{\mathrm{ad}}\triangleleft) \circ (_X\Pi \otimes \mathrm{id}_H) \,. \tag{20}$$

**Proof.** The first identity of (20) holds because $_X\Pi$ is an intertwining operator between the action $\mu_l$ and the trivial action. The second identity uses the Hopf module compatibility conditions and the fact that the antipode is a coalgebra anti-homomorphism. □

Hence we can define the right module structure on $_H X$ by

$$\mu_r^{_H X} := {}_X\mathrm{i} \circ ({}_{\mathrm{ad}}\triangleleft) \circ ({}_X\mathrm{p} \otimes \mathrm{id}_H) = {}_X\mathrm{i} \circ \mu_r \circ ({}_X\mathrm{p} \otimes \mathrm{id}_H). \tag{21}$$

In this way $_A(-)$ extends to a functor between $^H_H\mathcal{C}_H$ and $\mathcal{C}_H$. Similar to Theorem 3.5.2 one obtains

**Proposition 3.6.3** *Let $H$ be a Hopf algebra in $\mathcal{C}$. Then the categories $\mathcal{C}_H$ and $^H_H\mathcal{C}_H$ are monoidal equivalent. The equivalence is given by* $\mathcal{C}_H \xrightleftharpoons[_H(-)]{H\ltimes(-)} {^H_H\mathcal{C}_H}$.

# 4 Crossed Modules and Hopf Bimodules

In the first part of Section 4 we recall some results of [B2] on crossed modules in braided monoidal categories. Then Hopf bimodules in $\mathcal{C}$ will be introduced and a theorem will be proved which states the equivalence of crossed modules and Hopf bimodules. Crossed modules are (braided) analogues of the $k$-vector space of invariant one-forms of a bicovariant bimodule (Hopf bimodule) [Wor] and have a close connection to smash product constructions (see [Rad] in the case of bialgebras over the field $k$ and [B2] for the braided case).

## 4.1

**Definition 4.1.1** *[B2] Let $B$ be a bialgebra in the category $\mathcal{C}$. A right crossed module $(X, \mu_r, \nu_r)$ over $B$ is a right $B$-module and a right $B$-comodule obeying the compatibility relations*

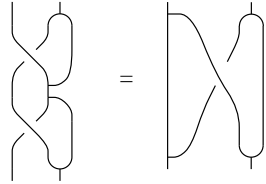

(Fig. 5)

*The category $\mathcal{DY}(\mathcal{C})^B_B$ is the category of crossed modules where the morphisms are both right module and right comodule morphisms over $B$. In a similar way all other combinations of crossed modules will be defined.*

**Example 4.1.1** *The unit object $(\underline{1}, \varepsilon_B, \eta_B)$ in $\mathcal{C}$ is a crossed module over $B$.*

**Example 4.1.2** *A Hopf algebra $H$ is an $H$-bicomodule through the regular coaction $\Delta$ and an $H$-bimodule through the regular action $\mathrm{m}$. Using the corresponding right adjoint action for $H$ according to (4) and the regular coaction, it follows that $(H, {}_{\mathrm{ad}}\triangleleft, \Delta)$ is a right crossed module over $H$. Dually analogue one can construct a crossed structure on $H$ by using the adjoint coaction and the regular action.*

The (pre-)braided monoidal structure of $\mathcal{DY}(\mathcal{C})^B_B$ is described in the following theorem.

**Theorem 4.1.1** *[B2] Let $B$ be a bialgebra in $\mathcal{C}$. Then the category $(\mathcal{DY}(\mathcal{C})^B_B, \otimes, \underline{1})$ is monoidal. It is pre-braided through*

$$^{\mathcal{DY}(\mathcal{C})^B_B}\Psi_{XY} := (\mathrm{id}_Y \otimes \mu_r^X) \circ (\Psi_{XY} \otimes \mathrm{id}_B) \circ (\mathrm{id}_X \otimes \nu_r^Y) \tag{22}$$

where $X, Y \in \mathrm{Ob}(\mathcal{DY}(\mathcal{C})_B^B)$. If $H$ is a Hopf algebra with isomorphic antipode in $\mathcal{C}$ then $\mathcal{DY}(\mathcal{C})_H^H$ is braided, i.e. the inverse of (22) exists and equals

$$(^{\mathcal{DY}(\mathcal{C})_B^B}\Psi_{XY})^{-1} = (\mu_r^X \otimes \mathrm{id}_Y) \circ (\mathrm{id}_Y \otimes (\Psi_{HY})^{-1}) \circ ((\Psi_{XY})^{-1} \otimes S^{-1}) \circ \quad (23)$$
$$\circ (\mathrm{id}_Y \otimes (\Psi_{XH})^{-1}) \circ (\nu_r^Y \otimes \mathrm{id}_X)$$

Consider the crossed modules over a Hopf algebra $H$ in $\mathcal{C}$ which has an isomorphic antipode. Then there are two canonical ways to construct left $H$-crossed modules from right $H$-crossed modules, and vice versa [B2]. For a right crossed module $(X, \mu_r, \nu_r)$ one obtains the left crossed modules $X^S = (X, (\mu^S)_l, (\nu^S)_l)$ and $^S X = (X, (^S\mu)_l, (^S\nu)_l)$ where the (co-)actions are defined through

$$\begin{aligned}
(\mu^S)_l &= \mu_r \circ \Psi_{XH}^{-1} \circ (S^{-1} \otimes \mathrm{id}_X), \\
(\nu^S)_l &= (S \otimes \mathrm{id}_X) \circ \Psi_{XH} \circ \nu_r, \\
\\
(^S\mu)_l &= \mu_r \circ \Psi_{HX} \circ (S \otimes \mathrm{id}_X), \\
(^S\nu)_L &= (S^{-1} \otimes \mathrm{id}_X) \circ \Psi_{HX}^{-1} \circ \nu_r.
\end{aligned} \quad (24)$$

Conversely for a left crossed module $(Y, \mu_l, \nu_l)$ one obtains the corresponding right crossed modules $Y_S = (Y, (\mu_S)_r, (\nu_S)_r)$ and $_S Y = (Y, (_S\mu)_r, (_S\nu)_r)$ by

$$\begin{aligned}
(\mu_S)_r &= \mu_l \circ \Psi_{YH} \circ (\mathrm{id}_Y \otimes S), \\
(\nu_S)_r &= (\mathrm{id}_Y \otimes S^{-1}) \circ \Psi_{YH}^{-1} \circ \nu_l, \\
\\
(_S\mu)_r &= \mu_l \circ \Psi_{HY}^{-1} \circ (\mathrm{id}_Y \otimes S^{-1}), \\
(_S\nu)_L &= (\mathrm{id}_Y \otimes S) \circ \Psi_{HX} \circ \nu_l.
\end{aligned} \quad (25)$$

## 4.2

In the next definition Hopf bimodules in the braided monoidal category $\mathcal{C}$ will be defined. In the case of bialgebras over a field $k$ they appear in the construction of bicovariant differential calculi over a Hopf algebra [Wor]. Through the generalization to braided categories (see below) it is possible to construct (bicovariant) differential calculi in such categories; this will be demonstrated in a forthcoming paper. See also [B3] on the construction of a first order differential calculus on the dual quantum braided matrix group.

**Definition 4.2.1** *Let $B$ be a bialgebra in $\mathcal{C}$. An object $(X, \mu_r, \mu_l, \nu_r, \nu_l)$ is called a $B$-Hopf bimodule if $(X, \mu_r, \mu_l)$ is a $B$-bimodule, and $(X, \nu_r, \nu_l)$ is a $B$-bicomodule in the category of $B$-bimodules, where the regular (co-)action on $B$ and the diagonal (co-)action on tensor products of modules are used. Hopf bimodules together with the $B$-bimodule-$B$-bicomodule morphisms form the category which will be denoted by $_B^B \mathcal{C}_B^B$.*

**Example 4.2.1** *Any bialgebra $B$ is a Hopf bimodule over itself if the regular actions and coactions $\mu_l := \mathrm{m} =: \mu_r$ and $\nu_l := \Delta =: \nu_r$ respectively are used.*

In the following proposition an analogue of Lemma 3.3.1 and Proposition 3.3.2 will be proved for crossed modules and Hopf bimodules.

**Proposition 4.2.1** *Let $X$ be a Hopf bimodule and $Y$ be a right crossed module over a bialgebra $B$ in $\mathcal{C}$. Then $X \otimes Y$ is a $B$-Hopf bimodule if it is equipped with the left action and coaction induced by $X$ and with the diagonal right action and coaction. In particular the cross product $B \ltimes Y$ [B2] is a Hopf bimodule. This induces two functors*

$$\otimes : {}_B^B \mathcal{C}_B^B \times \mathcal{DY}(\mathcal{C})_B^B \to {}_B^B \mathcal{C}_B^B \quad \text{and} \quad B \ltimes (-) : \mathcal{DY}(\mathcal{C})_B^B \to {}_B^B \mathcal{C}_B^B$$

*from which the functor $B \ltimes (-)$ is a full embedding.*

**Proof.** The bi-(co-)module axioms for $X\otimes Y$ are easily checked. By Lemma 3.3.1 $X\otimes Y$ is a left $B$-Hopf module. The compatibility conditions for the pairs $(\mu_l, \nu_r)$ and $(\mu_r, \nu_l)$ are straightforwardly verified by the use of the Hopf bimodule properties of $X$ and by (co-)associativity of $B$. The proof of the right $B$-Hopf module axiom is performed in Figure 7. In the first equation of Figure 7 it is made use of the Hopf bimodule properties of $X$ and the (co-)associativity of $B$. The result in the second equation can be shown with the help of the crossed module axiom for Y according to Figure 6. It is obvious that $\otimes$ and $B \ltimes (-)$ are functors. $B \ltimes (-)$ is a full embedding because of Proposition 3.3.2. $\square$

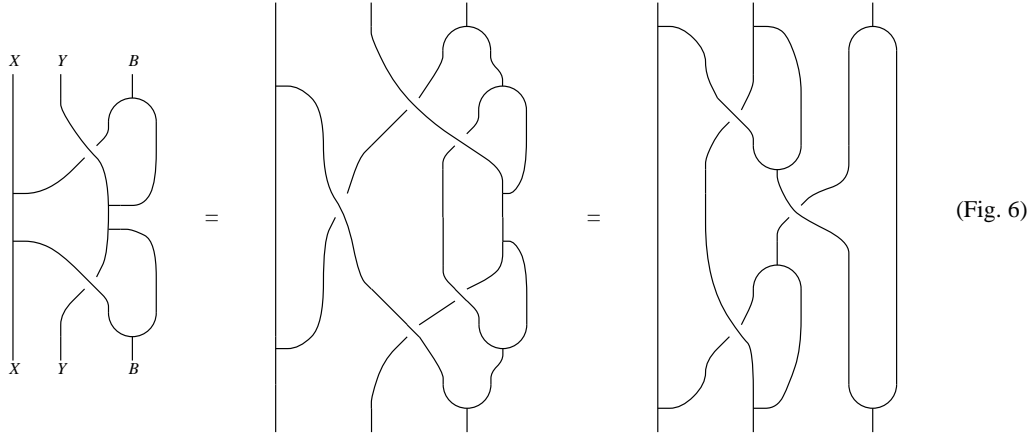

(Fig. 6)

**Proposition 4.2.2** *Let $(X, \mu_r, \mu_l, \nu_r, \nu_l)$ be a Hopf bimodule over a Hopf algebra $H$. Then the object $(X_{\mathrm{ad}}, {}_{\mathrm{ad}}\triangleleft, \nu_r)$ is a right $H$-crossed module, where the right adjoint action (4) is used. Dually analogue, $X$ together with its right action and the right coadjoint coaction is a right $H$-crossed module $(X^{\mathrm{ad}}, \mu_r, \nu_r^{\mathrm{ad}})$. The idempotent $_X\Pi : X^{\mathrm{ad}} \to X_{\mathrm{ad}}$ is crossed module morphism.*

**Proof.** The proof that $(X_{\mathrm{ad}}, {}_{\mathrm{ad}}\triangleleft, \nu_r)$ is a right $H$-crossed module is done in Figure 7 where in the second equation it is used that $X$ is a Hopf bimodule, $H$ is (co-)associative and the antipode is anti-comultiplicative. The antipode axioms and (co-)associativity of $H$ yield the result in the third equation.
The fact that $_X\Pi$ is crossed module morphism follows from (20) and its dual form. $\square$

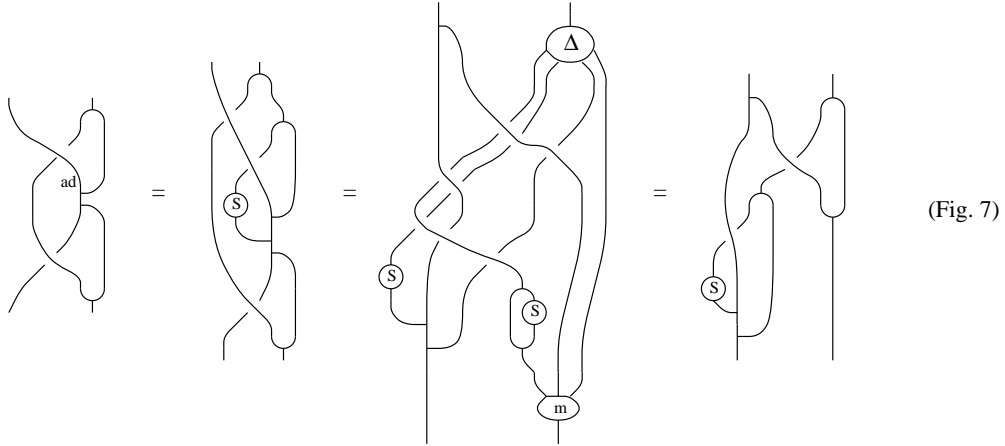

(Fig. 7)

**Remark 4.2.1** *The previous construction applied to the regular Hopf bimodule $H$ leads to the $H$-crossed modules $H_{\mathrm{ad}}$ and $H^{\mathrm{ad}}$. This yields new solutions of the Yang-Baxter equation for $H$ which arise through the braiding of these ojects in the category of crossed modules [Ma6].*

Proposition 4.2.2 will now be used to construct the functor $H \ltimes (-)$ from the Hopf bimodules into the crossed modules.

**Proposition 4.2.3** *Let $H$ be a Hopf algebra in $\mathcal{C}$. Then for an $H$-Hopf bimodule $X$ there exists a unique right $H$-crossed module structure on the object $_H X$ such that $_X\mathrm{p} : X^{\mathrm{ad}} \to {}_H X$ and $_X\mathrm{i}_H X \to X_{\mathrm{ad}}$ are crossed module morphisms. Explicitly the action is defined according to (21) and the coaction is defined in the dual symmetric manner.*
*As in the case of Hopf modules this construction induces a functor*

$$_H(-) : {}_H^H\mathcal{C}_H^H \longrightarrow \mathcal{DY}(\mathcal{C})_H^H$$

**Proof.** Using the equations (20) and Proposition 4.2.2 for the Hopf bimodule $X$ the statement follows directly. □

### 4.3

The monoidal structure of the Hopf bimodules can be directly deduced from Theorem 3.5.1. For the braiding however one has to take into account the compatibility of the right (co-)module structures.

**Theorem 4.3.1** *Let $H$ be a Hopf algebra in the category $\mathcal{C}$. Then the category of Hopf bimodules over $H$ is monoidal with the tensor product given by*

$$\otimes_H := \otimes \circ (\mathrm{id} \times {}_H(-)) : {}_H^H\mathcal{C}_H^H \times {}_H^H\mathcal{C}_H^H \longrightarrow {}_H^H\mathcal{C}_H^H,$$

*with the unit object $H$ as regular $H$-Hopf bimodule, and with the natural isomorphisms ${}^{{}_H^H\mathcal{C}_H^H}\alpha$, ${}^{{}_H^H\mathcal{C}_H^H}\rho$ and ${}^{{}_H^H\mathcal{C}_H^H}\lambda$ according to*

$$({}^{{}_H^H\mathcal{C}_H^H}\alpha)_{X\,Y\,Z} := \alpha_{X\,{}_H Y\,{}_H Z}\,, \quad ({}^{{}_H^H\mathcal{C}_H^H}\rho)_X := \rho_X\,, \quad ({}^{{}_H^H\mathcal{C}_H^H}\lambda)_X := {}_X\mu$$

*where $X, Y, Z \in \mathrm{Ob}({}_H^H\mathcal{C}_H^H)$ and $_X\mu$ is given by (8). ${}_H^H\mathcal{C}_H^H$ is pre-braided with the pre-braiding given on the objects $X, Y \in \mathrm{Ob}({}_H^H\mathcal{C}_H^H)$ through*

$${}^{{}_H^H\mathcal{C}_H^H}\Psi_{X\,Y} = (\mu_l^Y \otimes {}_X\mathrm{p} \circ \mu_r^X) \circ (\mathrm{id}_H \otimes \Psi_{X\,Y} \otimes \mathrm{id}_H) \circ (\nu_l^X \otimes \nu_r^Y \circ {}_Y\mathrm{i}) \,. \tag{26}$$

*The category ${}_H^H\mathcal{C}_H^H$ is braided if the antipode of $H$ is an isomorphism.*

**Proof.** The functor $\otimes_H$ is well defined by Propositions 4.2.1 and 4.2.3. $H$ is an $H$-Hopf bimodule. To prove the monoidal structure of ${}_H^H\mathcal{C}_H^H$ it is therefore sufficient by Theorem 3.5.1 to prove the right (co-)module property of the morphisms $\alpha$, $\rho$ and $\lambda$. This is a straightforward calculation.

It can be verified explicitly that ${}^{{}_H^H\mathcal{C}_H^H}\Psi$ given by (26) is a (pre-)braiding. However anticipating the result of Theorem 4.3.2 on the monoidal equivalence of ${}_H^H\mathcal{C}_H^H$ and $\mathcal{DY}(\mathcal{C})_H^H$ – for its proof only the momoidal structures of both categories are needed – one obtains by Theorem 4.1.1 that for ${}^{{}_H^H\mathcal{C}_H^H}\Psi$ in eq. (26) it holds

$${}^{{}_H^H\mathcal{C}_H^H}\Psi_{X\,Y} = ({}_Y\mu \otimes {}_X\mu) \circ (\xi_{{}_H Y\,{}_H X})^{-1} \circ (H \ltimes {}^{\mathcal{DY}(\mathcal{C})_H^H}\Psi_{{}_H X\,{}_H Y}) \circ \xi_{{}_H X\,{}_H Y} \circ ({}_X\nu \otimes {}_Y\nu) \tag{27}$$

and hence ${}^{{}_H^H\mathcal{C}_H^H}\Psi$ is a braiding in ${}_H^H\mathcal{C}_H^H$. Here $_X\nu$ and $_X\mu$ of Lemma 3.3.3 have been used which apply to Hopf bimodules and induce natural isomorphisms. The natural isomorphism $\xi$ is defined in the proof of Theorem 3.5.2. □

The next theorem shows that the braided monoidal structures of ${}_H^H\mathcal{C}_H^H$ and $\mathcal{DY}(\mathcal{C})_H^H$ are compatible with the functors $_H(-)$ and $H \ltimes (-)$.

**Theorem 4.3.2** *Let $H$ be a Hopf algebra in $\mathcal{C}$ with isomorphic antipode. Then the categories $\mathcal{DY}(\mathcal{C})_H^H$ and ${}_H^H\mathcal{C}_H^H$ are braided monoidal equivalent. The equivalence is given by $\mathcal{DY}(\mathcal{C})_H^H \xrightleftharpoons[{}_H(-)]{H\ltimes(-)} {}_H^H\mathcal{C}_H^H$.*

**Proof.** To prove the monoidal equivalence of $\mathcal{DY}(\mathcal{C})_H^H$ and ${}_H^H\mathcal{C}_H^H$ it is sufficient to prove that the natural isomorphisms $\lambda$ and $\mu$, which appear in the proof of Theorem 3.5.2, are natural isomorphisms if they are restricted to $\mathcal{DY}(\mathcal{C})_H^H$ and ${}_H^H\mathcal{C}_H^H$. In other words it has to be verified that $\lambda_X$ and ${}_X\mu$ are right (co-)module morphisms which is a direct and simple calculation.

The compatibility of the braidings ${}^{{}_H^H\mathcal{C}_H^H}\Psi$ and ${}^{\mathcal{DY}(\mathcal{C})_B^B}\Psi$ with the monoidal equivalent structures of $\mathcal{DY}(\mathcal{C})_H^H$ and ${}_H^H\mathcal{C}_H^H$ holds by construction of ${}^{{}_H^H\mathcal{C}_H^H}\Psi$ (see the argumentation in the proof of Theorem 4.3.1). $\square$

**Remark 4.3.1** *The mirror symmetric construction (using left crossed modules) leads to another braiding on ${}_H^H\mathcal{C}_H^H$ which is given by a diagram similar to (26) but with the replacement $\mu \leftrightarrow \mathrm{p} \circ \mu$ and $\nu \leftrightarrow \nu \circ \mathrm{i}$. This braiding is equivalent to the braiding in (26).*

**Remark 4.3.2** *To get in touch with the results of [Sch] for the purely symmetric case, one observes that the braided version of the morphism in [Sch] which rules the braiding of $H$-Hopf bimodules (in the symmetric case) can be written as*

$$\Psi'_{XY} = (\mu_l^Y \otimes \mu_r^X) \circ (\mathrm{id}_H \otimes (\Psi_{XY} \circ ({}_X\Pi \otimes \Pi_Y)) \otimes \mathrm{id}_H) \circ (\nu_l^X \otimes \nu_r^Y). \tag{28}$$

*Using Proposition 3.4.1 one finds that for $\lambda_{Y,X}^H \circ \Psi'_{XY}$ the factorization $\lambda_{Y,X}^H \circ \Psi'_{XY} = {}^{{}_H^H\mathcal{C}_H^H}\Psi_{XY} \circ \lambda_{X,Y}^H$ over $\lambda_{X,Y}^H$ holds which implies that $\lambda_{Y,X}^H \circ \Psi'_{XY}$ coequalizes the morphisms $\mu_r^X \otimes \mathrm{id}_Y$ and $\mathrm{id}_X \otimes \mu_l^Y$ and that the factorization is unique.*

# 5 Bialgebra Projections and Cross Products

In this section the connection of bialgebra projections and Hopf bimodule bialgebras is investigated. The resulting theorem is applied to bialgebra cross products in the braided monoidal category $\mathcal{C}$.

## 5.1

**Definition 5.1.1** *Let $H$ be a Hopf algebra in $\mathcal{C}$ and $(X, \mu_r, \mu_l, \nu_r, \nu_l)$ be an $H$-Hopf bimodule. Then the relative antipode $S_{X/H}$ of $X$ w.r.t. $H$ is defined by*

$$S_{X/H} := M_X \circ (S \otimes \mathrm{id}_X \otimes S) \circ N_X \tag{29}$$

*where $M_X := \mu_l \circ (\mathrm{id}_H \otimes \mu_r)$ and $N_X := (\mathrm{id}_H \otimes \nu_r) \circ \nu_l$.*

**Proposition 5.1.1** *Let $H$ and $X$ be like in Definition 5.1.1.*

1. *The "polarized" form of the anti-(co-)multiplicity of the antipode holds.*

   $$S_{X/H} \circ \mu_r = \mu_l \circ \Psi_{XH} \circ (S_{X/H} \otimes S), \quad S_{X/H} \circ \mu_l = \mu_r \circ \Psi_{HX} \circ (S \otimes S_{X/H}).$$

   *And dually analogue for the coactions.*

2. *The relative antipode $S_{H/H}$ coincides with the antipode $S$.*

3.

$$S_{X/H} \circ \Pi_X = {}_X\Pi \circ \Pi_X = {}_X\Pi \circ S_{X/H}$$
$$S_{X/H} \circ {}_X\Pi = \Pi_X \circ {}_X\Pi = \Pi_X \circ S_{X/H}$$

4. *Suppose that the antipode $S$ of $H$ is isomorphic. Then the relative antipode $S_{X/H}$ is an isomorphism. Its inverse is given by*

$$S_{X/H}^{-1} = M_X \circ (S^{-1} \otimes (\Psi_{X\,H})^{-1}) \circ ((\Psi_{H\,H})^{-1} \otimes \mathrm{id}_X) \circ (S^{-1} \otimes (\Psi_{H\,X})^{-1}) \circ N_X$$

where $M_X$ and $N_X$ of Definition 5.1.1 are used.

5. *The morphisms ${}_X S := \mathrm{p}_X \circ S_{X/H} \circ {}_X\mathrm{i}$ and $S_X := {}_X\mathrm{p} \circ S_{X/H} \circ \mathrm{i}_X$ are module morphisms ${}_X S : {}_H X \to (X_H)_S$ and $S_X : X_H \to {}^S({}_H X)$, and comodule morphisms ${}_X S : {}_H X \to {}_S(X_H)$ and $S_X : X_H \to ({}_H X)^S$ according to (24) and (25).*

**Proof.** The Hopf bimodule properties, Proposition 3.2.1 (and its mirror reversed version) and the deduced identities $({}_X\Pi \otimes \mathrm{id}_H) \circ \nu_r^X \circ {}_X\mathrm{i} = \nu_r^X \circ {}_X\mathrm{i}$ and $\mathrm{p}_X \circ \mu_l^X \circ (\mathrm{id}_H \otimes \Pi_X) = \mathrm{p}_X \circ \mu_l^X$ immediately imply statements (1), (2), (3) and (5). Statement (1) and the equality $\mu_l^X \circ (\mathrm{id}_H \otimes S_{X/H}) \circ \nu_l^X = \Pi_X$ yield statement (4). □

**Proposition 5.1.2** *Let $H$ be a Hopf algebra in $\mathcal{C}$ and let $X$ and $Y$ be $H$-Hopf bimodules. Then the following identity and its dual symmetric version hold.*

$$ {}_A^A \mathcal{C}_A^A \Psi_{X,Y} \circ S_{X \otimes_A Y/A} \circ \underset{A}{\otimes} = \underset{A}{\otimes} \circ {}^{\mathcal{C}} \Psi_{X,Y} \circ (S_{X/A} \otimes S_{Y/A}). \tag{30}$$

In the next definition bialgebra projections will be considered. They turn out to be the key ingredient for cross product constructions. It will be shown that Hopf bimodule bialgebras are essentially bialgebra projections.

**Definition 5.1.2** *Let $B$ and $H$ be bialgebras in $\mathcal{C}$ and $H \xrightarrow{\underline{\eta}} B \xrightarrow{\underline{\varepsilon}} H$ be a pair of bialgebra morphisms such that $\underline{\varepsilon} \circ \underline{\eta} = \mathrm{id}_H$. Then $(H, B, \underline{\eta}, \underline{\varepsilon})$ is called a bialgebra projection on $H$.*

It is not difficult to verify

**Lemma 5.1.3** *Let $(H, B, \underline{\eta}, \underline{\varepsilon})$ be a bialgebra projection on $H$ in $\mathcal{C}$. Then $\underline{B} = (B, \mu_r^B, \mu_l^B, \nu_r^B, \nu_l^B)$ is an $H$-Hopf bimodule where*

$$\mu_l^B = \mathrm{m}_B \circ (\underline{\eta} \otimes \mathrm{id}_B), \quad \mu_r^B = \mathrm{m}_B \circ (\mathrm{id}_B \otimes \underline{\eta}), \tag{31}$$
$$\nu_l^B = (\underline{\varepsilon} \otimes \mathrm{id}_B) \circ \Delta_B, \quad \nu_r^B = (\mathrm{id}_B \otimes \underline{\varepsilon}) \circ \Delta_B.$$

The following proposition recovers the deep connection of Hopf bimodule bialgebras and bialgebra projections.

**Proposition 5.1.4** *Let $H$ be a Hopf algebra in $\mathcal{C}$. Then for every bialgebra projection $(H, B, \underline{\eta}_B, \underline{\varepsilon}_B)$ the object $\underline{B}$ according to Lemma 5.1.3 is a bialgebra $(\underline{B}, \underline{\mathrm{m}}_B, \underline{\eta}_B, \underline{\Delta}_B, \underline{\varepsilon}_B)$ in ${}_H^H \mathcal{C}_H^H$ where the projection morphisms $\underline{\eta}_B$ and $\underline{\varepsilon}_B$ are the unit and counit respectively, and with the multiplication $\underline{\mathrm{m}}_B$ and the comultiplication $\underline{\Delta}_B$ defined through*

$$\underline{\mathrm{m}}_B = \mathrm{m}_B \circ (\mathrm{id}_B \otimes {}_H\mathrm{i}) \quad and \quad \underline{\Delta}_B = (\mathrm{id}_B \otimes {}_H\mathrm{p}) \circ \Delta. \tag{32}$$

*Conversely every bialgebra $\underline{B} = (B, \underline{\mathrm{m}}_B, \underline{\eta}_B, \underline{\Delta}_B, \underline{\varepsilon}_B)$ in ${}_H^H \mathcal{C}_H^H$ can be turned into a bialgebra $B = (B, \mathrm{m}_B, \eta_B, \Delta_B, \varepsilon_B)$ in $\mathcal{C}$ where the structure morphisms are given by*

$$\mathrm{m}_B = \underline{\mathrm{m}}_B \circ \lambda_{B,B}^H, \quad \eta_B = \underline{\eta}_B \circ \eta_H, \quad \Delta_B = \rho_{B,B}^H \circ \underline{\Delta}_B, \quad \varepsilon_B = \varepsilon_H \circ \underline{\varepsilon}_B. \tag{33}$$

*Furthermore* $(H, B, \underline{\eta}_B, \underline{\varepsilon}_B)$ *is a bialgebra projection on* $H$.
*There also exists a correspondence of Hopf algebra structures. For any Hopf algebra projection* $(H, B, \underline{\eta}_B, \underline{\varepsilon}_B)$ *on* $H$ *the antipode of* $\underline{B}$ *is given by*

$$\underline{S}_B = M_B \circ (\mathrm{id}_H \otimes S_B \otimes \mathrm{id}_H) \circ N_B \tag{34}$$

*and for any Hopf algebra* $\underline{B}$ *in* $^H_H\mathcal{C}^H_H$ *the antipode of* $B$ *is given by*

$$S_B = \underline{S}_B \circ S_{B/H} = S_{B/H} \circ \underline{S}_B \,. \tag{35}$$

**Proof.** Suppose $(H, B, \underline{\eta}_B, \underline{\varepsilon}_B)$ is a bialgebra projection on $H$. Then by Lemma 5.1.3 $\underline{B}$ is an $H$-Hopf bimodule. The proof that $\underline{m}_B, \underline{\eta}_B, \underline{\Delta}_B, \underline{\varepsilon}_B$ and $\underline{S}_B$ are Hopf bimodule morphisms essentially makes use of Proposition 3.2.1 and the fact that $\underline{\eta}_B$ and $\underline{\varepsilon}_B$ are bialgebra morphisms in $\mathcal{C}$ which fulfill $\underline{\varepsilon}_B \circ \underline{\eta}_B = \mathrm{id}_H$. The algebra and coalgebra axioms can be verified straigthforwardly. In case $B$ is a Hopf algebra the proof of the antipode axioms (of $\underline{B}$) uses the bi-(co-)module properties of $\underline{B}$ and the Hopf axioms of $H$ and $B$. The proof of the bialgebra axiom

$$\underline{\Delta}_B \circ \underline{m}_B = (\underline{m}_B \otimes_H \underline{m}_B) \circ (\mathrm{id}_B \otimes_H {}^{A}_{A}\mathcal{C}^A_A \Psi_{B,B} \otimes_H \mathrm{id}_B) \circ (\underline{\Delta}_B \otimes_H \underline{\Delta}_B) \tag{36}$$

needs the identity

$$(\nu_r^B \circ {}_B\Pi \otimes \mathrm{id}_B) \circ \Delta_B \circ {}_B\mathrm{i} = (\mathrm{id}_B \otimes \nu_l^B) \circ \Delta_B \circ {}_B\mathrm{i} \tag{37}$$

and is indicated in Figure 8 where the morphisms ${}_B\Pi$, ${}_B\mathrm{i}$ and ${}_B\mathrm{p}$ are used. The left hand side of (36) is transformed to the left diagram of Figure 8 by use of (26) and Proposition 3.2.1. In the first equation of Figure 8 the equality (37) is exploited. To verify the second equation again Proposition 3.2.1 is taken into account.

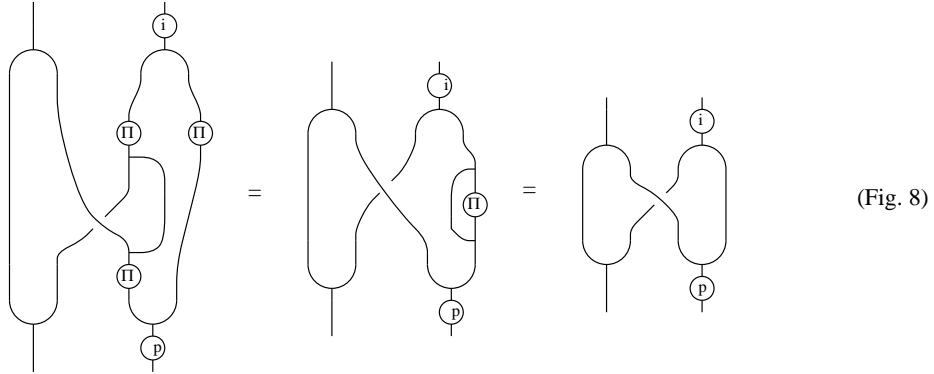

(Fig. 8)

Conversely let $\underline{B} = (B, \underline{m}_B, \underline{\eta}_B, \underline{\Delta}_B, \underline{\varepsilon}_B\,(, \underline{S}_B))$ be a bialgebra (Hopf algebra) in $^H_H\mathcal{C}^H_H$. The verification that $(H, B, \underline{\eta}_B, \underline{\varepsilon}_B)$ is a bialgebra projection on $H$ is simple. It is straightforward to prove that $(B, \mathrm{m}_B, \eta_B)$ and $(B, \Delta_B, \varepsilon_B)$ according to (33) are algebra and coalgebra respectively. A little more subtle is the proof of the Hopf axioms for the antipode (35). The most difficult part is the verification of the multiplicativity of the comultiplication $\Delta_B$ in $\mathcal{C}$. The corresponding proof is sketched in Figure 9 where the multiplication $\underline{m}_B$, the comultiplication $\underline{\Delta}_B$ and the morphisms ${}_B\Pi$, ${}_B\mathrm{i}$ and ${}_B\mathrm{p}$ are used in particular. With the help of Proposition 3.4.1 and the bialgebra properties of $\underline{B}$ a lengthy calculation transforms $\Delta_B \circ \mathrm{m}_B$ to the left diagram in Figure 9. In the first equation of Figure 9 it is used that $\underline{B}$ is a Hopf bimodule. The second equation exploits again Hopf bimodule properties and Proposition 3.2.1. □

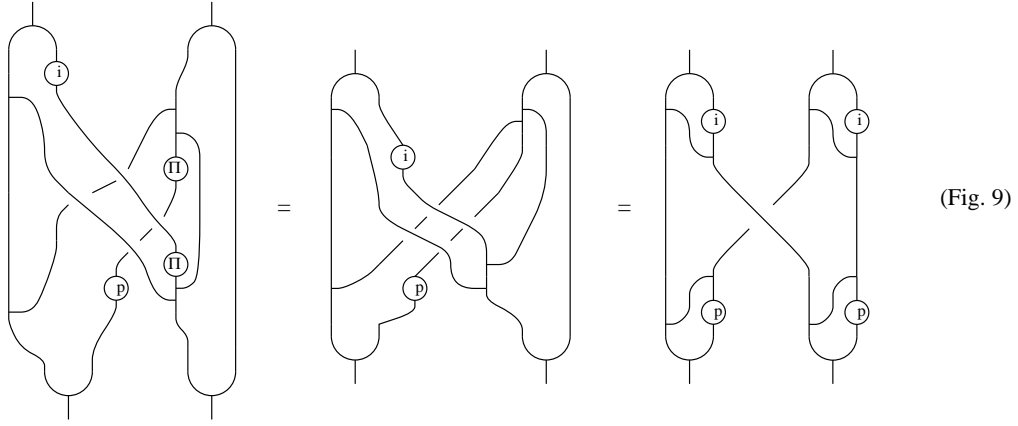
(Fig. 9)

For a Hopf algebra $H$ in $\mathcal{C}$ we denote by $H$-Bialg-$\mathcal{C}$ / $H$-Hopf-$\mathcal{C}$ the category of bialgebra / Hopf algebra projections on $H$. Its objects are the projections $B = \{B, \mu_B, \eta_B, \Delta_B, \varepsilon_B, (S_B), \underline{\eta}_B, \underline{\varepsilon}_B\}$ and its morphisms are bialgebra / Hopf algebra morphisms $f : B \to D$ such that $f \circ \underline{\eta}_B = \underline{\eta}_D$ and $\underline{\varepsilon}_D \circ f = \underline{\varepsilon}_B$. Then the previous proposition can be used to formulate

**Theorem 5.1.5** *The categories $H$-Bialg-$\mathcal{C}$ and Bialg-${}_H^H\mathcal{C}_H^H$ are isomorphic. The isomorphism is given by* Bialg-${}_H^H\mathcal{C}_H^H \xrightarrow[F]{G} H$-Bialg-$\mathcal{C}$ *where the functors $F$ and $G$ are defined through Proposition 5.1.4 on the objects and through the identity on the morphisms.*

**Proof.** It is not difficult to prove that the functors $F$ and $G$ are well defined. After some calculation the constructions given in Proposition 5.1.4 for the objects turn out to be inverse to each other, and since $F$ and $G$ are identities on the morphisms this proves the theorem. □

## 5.2

The equivalence of the braided categories ${}_H^H\mathcal{C}_H^H$ and $\mathcal{DY}(\mathcal{C})_H^H$ allows us to reformulate the previous theorem in terms of the category $\mathcal{DY}(\mathcal{C})_H^H$. Then we immediately get, in one direction, the formulas of [B2] for cross products of Hopf algebras (bialgebras), and in the other direction the Radford-Majid criterion when a Hopf algebra is a cross product (see [Rad, Ma3], and [B2] in the braided case). This seems to be the most natural and easy description of the subject. For completeness it will be outlined in more detail in the following.

Very similar to the well known case of Hopf algebras over a field $k$ the smash product and the smash coproduct in the category $\mathcal{C}$ can be defined in a formal manner as a certain (co-)limit. For a Hopf algebra $H$ in $\mathcal{C}$ and an algebra $(A, \mu_A)$ in $\mathcal{C}_H$ the smash product $H \ltimes_{\mu_A} A$ is defined as the universal algebra in $\mathcal{C}$ such that:

1. There are algebra morphisms j : $H \to H \ltimes_{\mu_A} A$ and i : $A \to H \ltimes_{\mu_A} A$. In addition i is algebra morphism in $\mathcal{C}_H$ where the module structure on $H \ltimes_{\mu_A} A$ is the right adjoint action induced by the morphism j, i.e. i $\in \mathrm{Alg}_H(A, (H \ltimes_{\mu_A} A)_j)$

2. If $U$ is any algebra in $\mathcal{C}$ and $g \in \mathrm{Alg}(H, U)$, $f \in \mathrm{Alg}_H(A, U_g)$, then there exists a unique algebra morphism $g \ltimes_{\mu_A} f : H \ltimes_{\mu_A} A \to B$ such that $f = (g \ltimes_{\mu_A} f) \circ \mathrm{i}$ and $g = (g \ltimes_{\mu_A} f) \circ \mathrm{j}$.

The smash coproduct $H \ltimes^{\nu_C} C$ is defined in the obvious dual symmetric manner. Both products are unique up to (co-)algebra isomorphism. The smash product can be realized on the tensor product $H \otimes A$ through [Ma4, Dra]

$$\mathrm{m}_\ltimes = (\mathrm{m}_H \otimes \mathrm{m}_A) \circ (\mathrm{id}_H \otimes (\mathrm{id}_H \otimes \mu_r^A) \circ (\Psi_{AH} \otimes \mathrm{id}_H) \circ (\mathrm{id}_A \otimes \Delta_H) \otimes \mathrm{id}_A),$$
$$\eta_\ltimes = \eta_H \otimes \eta_A, \quad \mathrm{i} = \eta_H \otimes \mathrm{id}_A, \quad \mathrm{j} = \mathrm{id}_H \otimes \eta_A$$

and $g \ltimes_{\mu_A} f = \mathrm{m}_U \circ (g \otimes f)$ for the corresponding unique morphism.

If $X$ is at the same time a right $H$-module algebra and a right $H$-comodule coalgebra such that the smash product and the smash coproduct (realized on $H \otimes X$) are compatible in such a way that $H \ltimes X := (H \otimes X, \mathrm{m}_\ltimes, \eta_\ltimes, \Delta_\ltimes, \varepsilon_\ltimes)$ is a bialgebra in $\mathcal{C}$ then we say, in the sense of [Rad], that the pair $(H, X)$ is *admissible*; $X$ is called $H$-*admissible object* in $\mathcal{C}$. The category $H$-cp-$\mathcal{C}$ is the category of admissible pairs $(H, X)$ with bialgebra morphisms $h : H \ltimes X \to H \ltimes Y$ such that $h \circ \mathrm{j}_X = \mathrm{j}_Y$ and $\mathrm{k}_Y \circ h = \mathrm{k}_X$, where $\mathrm{k} = \varepsilon \otimes \mathrm{id}_H$.

Like in [Rad] one proves that the following relations hold for an admissible pair $(H, X)$ in $\mathcal{C}$.

$$\begin{aligned}
\varepsilon_X \circ \mathrm{m}_X &= \varepsilon_X \otimes \varepsilon_X, & \varepsilon_X \circ \mu_r^X &= \varepsilon_X \otimes \varepsilon_H, \\
\Delta_X \circ \eta_X &= \eta_X \otimes \eta_X, & \nu_r^X \circ \eta_X &= \eta_X \otimes \eta_H, \\
\varepsilon_X \circ \eta_X &= \mathrm{id}_{\underline{1}}
\end{aligned} \tag{38}$$

From (38) it follows that $(H, H \ltimes X, \mathrm{id}_H \otimes \eta_X, \mathrm{id}_H \otimes \varepsilon_X)$ is a bialgebra projection on $H$. Using the results of [B2] and Theorems 4.3.2 and 5.1.5 this leads to the description of $H$-admissible objects in the category $\mathcal{C}$ in terms of crossed module bialgebras.

**Theorem 5.2.1** *Let $H$ be a Hopf algebra in $\mathcal{C}$ with isomorphic antipode. Then the category of admissible pairs $H$-cp-$\mathcal{C}$ and the category of $H$-crossed module bialgebras* Bialg-$\mathcal{DY}(\mathcal{C})_H^H$ *are isomorphic. The isomorphism is given through* Bialg-$\mathcal{DY}(\mathcal{C})_H^H \xrightarrow[{}_H(-)\circ F]{H \ltimes (-)} H$-cp-$\mathcal{C}$.

**Proof.** Let $X \in \mathrm{Ob}(\mathrm{Bialg}\text{-}\mathcal{DY}(\mathcal{C})_H^H)$ then $H \ltimes X \in \mathrm{Ob}(H\text{-cp-}\mathcal{C})$ according to [B2]. Using Theorem 5.1.5 and Proposition 4.2.3 yields ${}_H F(H \ltimes X) = X$ as bialgebra in $\mathcal{DY}(\mathcal{C})_H^H$. Similarly ${}_H F(H \ltimes f) = f$ for $H$-crossed module bialgebra morphisms. Conversely let $H \ltimes Y \in \mathrm{Ob}(H\text{-cp-}\mathcal{C})$. Then according to the statement following (38) and Theorem 5.1.5 ${}_H F(H \ltimes Y) = Y$ is a bialgebra in $\mathcal{DY}(\mathcal{C})_H^H$ and hence $H \ltimes ({}_H F(H \ltimes Y)) = H \ltimes Y$ is in $H$-cp-$\mathcal{C}$ [B2]. Using that ${}_H \mathrm{p} = \varepsilon_H$, ${}_H \mathrm{i} = \eta_H$ and the fact that every crossed module morphism $f : H \ltimes Y \to H \ltimes Z$ factorizes (see the proof of Proposition 3.3.2) yields $H \ltimes ({}_H F(H \ltimes f)) = H \ltimes f$. □

**Remark 5.2.1** *Analogous results like in Theorem 5.2.1 are obtained in the Hopf algebra case for Hopf-admissible pairs and crossed module Hopf algebras.*

**Remark 5.2.2** *Besides Theorem 5.2.1 another application of Theorem 5.1.5 is the construction of tensor algebras in $\mathcal{C}$ related to Hopf module tensor products in ${}_H^H \mathcal{C}_H^H$. This will be exploited in a forthcoming paper to construct (bicovariant) differential calculi in braided tensor categories.*

# A  Canonical splitting of idempotents

For any category $\mathcal{C}$ we consider the category $\widehat{\mathcal{C}}$ according to the definition in Section 3.

**Definition A.0.1** *We say that idempotents canonically split in a (braided) monoidal category $\mathcal{C}$ if the following data are given: a (braided) monoidal functor $\widehat{\mathcal{C}} \xrightarrow{s} \mathcal{C}$ such that $s \circ i = \mathrm{id}_\mathcal{C}$, and an isomorphism of the (braided) monoidal functors $i \circ s$ and $\mathrm{id}_{\widehat{\mathcal{C}}}$.*

**Lemma A.0.2** *There exists a natural canonical splitting of idempotents in $\widehat{\mathcal{C}}$. Any idempotent $f \in \mathrm{End}_{\widehat{\mathcal{C}}}(X_e)$ is also an idempotent in $\mathrm{End}_\mathcal{C}(X)$, the "splitting" functor $s : \widehat{\widehat{\mathcal{C}}} \to \widehat{\mathcal{C}}$ assigns to an object $(X_e)_f$ of $\widehat{\widehat{\mathcal{C}}}$ the object $X_f$ of $\widehat{\mathcal{C}}$ and identifies any set $\widehat{\widehat{\mathcal{C}}}((X_{e_1})_{f_1}, (X_{e_2})_{f_2})$ with the set $\widehat{\mathcal{C}}(X_{f_1}, X_{f_2})$ in a natural way. The functor transformation $i \circ s \to \mathrm{id}_{\widehat{\widehat{\mathcal{C}}}}$ is defined by the collection of isomorphisms $f : (X_e)_f \to (X_e)_f$ for any $(X_e)_f$ in $\widehat{\widehat{\mathcal{C}}}$.*

**Example A.0.1** *There exist two canonical splittings of idempotents for the category* $\mathrm{Vect}_k$ *of vector spaces over a field $k$ (as well as for any abelian category). With any idempotent $e \in \mathrm{End}_k(X)$ one connects two isomorphic objects $\mathrm{im}\, e$ (a subspace of $X$) and $\mathrm{coim}\, e$ (a factorspace of $X$). Any linear transformation $g \in \mathrm{Hom}_k(X,Y)$, with $fge = g$, where $f \in \mathrm{End}_k(Y)$ is another idempotent, defines the corresponding linear transformations $\mathrm{im}\, e \to \mathrm{im}\, f$ and $\mathrm{coim}\, e \to \mathrm{coim}\, f$.*

**Lemma A.0.3** *Let $\mathcal{C}$ be a monoidal category with canonical splitting of idempotents and $A$ be an algebra in $\mathcal{C}$. Then one can define the corresponding canonical splitting of idempotents on the category ${}_A\mathcal{C}$ of $A$-modules. For the categories of comodules, Yetter-Drinfel'd modules, and Hopf (bi-)modules the construction is similar.*

National Academy of Sciences
Bogolyubov Institute for Theoretical Physics
252 143, Kiev-143
Ukraine
e-mail: *mmtpitp@gluk.apc.org*

Department of Mathematics and Computer Science
University of Amsterdam
Plantage Muidergracht 24
NL-1018 TV Amsterdam
The Netherlands
e-mail: *drabant@fwi.uva.nl*